\documentclass[11pt, leqno]{article}
\usepackage{}
\usepackage[total={6in, 9in}]{geometry}
 \usepackage{amsfonts}
\usepackage{amsthm}
\usepackage{amssymb}
\usepackage{amsmath,amsfonts}
\usepackage{mathrsfs}
\usepackage{cases}
\usepackage{latexsym,bm}
\usepackage{indentfirst}
\usepackage{color}
\usepackage{ifpdf}
\usepackage{graphicx}
\usepackage{psfrag}
\usepackage{dsfont}
\usepackage[ruled,vlined]{algorithm2e}
\usepackage{enumerate}
\numberwithin{equation}{section}
\usepackage[
pdfauthor={CESS},
pdftitle={Spanning trails},
pdfstartview=XYZ,
bookmarks=true,
colorlinks=true,
linkcolor=blue,
urlcolor=blue,
citecolor=blue,
bookmarks=true,
linktocpage=true,
hyperindex=true
]{hyperref}

\usepackage{pgf,tikz,pgfplots}
\usepackage{tikz-3dplot}
\usepackage{pifont}
\usepackage{refcount}

\newtheorem{THM}{\textbf{Theorem}}[section]

\newtheorem{LEM}[THM]{\textbf{Lemma}}

\newtheorem{CON}[THM]{\textbf{Conjecture}}

\newcommand{\pf}{\noindent\textbf{Proof}.\quad}

\newtheorem*{THM1}{\textbf{Theorem 1.2}}
\newtheorem*{THM2}{\textbf{Theorem 1.3}}

\newcommand{\ve}{\varepsilon }

\newcommand{\CC}{\mathcal{C}}

\DeclareMathOperator{\df}{def}
\newcommand{\pbar}{\overline{\varphi}}

\newcommand{\arxiv}[1]{\href{http://arxiv.org/abs/#1}{\texttt{arXiv:#1}}}

\begin{document}
\title{Edge coloring graphs with large minimum degree}
\author{Michael J. Plantholt and Songling Shan \\ 
	\medskip  Illinois State  University, Normal, IL 61790\\
	\medskip 
	{\tt mikep@ilstu.edu; sshan12@ilstu.edu}
}

\date{\today}
\maketitle

\emph{\textbf{Abstract}.}
Let $G$ be a simple graph with maximum degree $\Delta(G)$. A subgraph $H$ of $G$ is overfull if $|E(H)|>\Delta(G)\lfloor |V(H)|/2 \rfloor$. Chetwynd and Hilton in 1985 conjectured that a graph $G$  with $\Delta(G)>|V(G)|/3$
has chromatic index $\Delta(G)$ if and only if $G$ contains no overfull subgraph. The 1-factorization conjecture is a special case of this overfull conjecture,
which states that for even $n$, every regular $n$-vertex graph with degree at least about $n/2$ has a 1-factorization and was confirmed for large graphs in 2014.  Supporting the overfull conjecture as well as generalizing the 1-factorization conjecture in an asymptotic  way, 
in this paper, we show that for any given $0<\ve <1$, there exists a positive integer $n_0$
such that the following statement holds:
if $G$ is a graph on $2n\ge n_0$ vertices with minimum degree at least $(1+\ve)n$,   then $G$ has chromatic index $\Delta(G)$ if and only if $G$ contains no overfull subgraph.


\vspace{2mm}

\section{Introduction}

In this paper, a graph means a simple graph
and a multigraph may contain parallel edges but no loops. 
Let $G$ be a multigraph.
Denote by $V(G)$ and  $E(G)$ the vertex set and edge set of $G$,
respectively, and by $e(G)$ the cardinality of $E(G)$. 
For $v\in V(G)$, $N_G(v)$ is the set of neighbors of $v$ 
in $G$, and 
$d_G(v)$, the degree of $v$
in $G$, is the number of edges of $G$ that are incident with $v$.
When $G$ is simple, $d_G(v)=|N_G(v)|$.  
For 
$S\subseteq V(G)$, $N_G(S)=\bigcup_{v\in S}N_G(v)$,  
 the subgraph of $G$ induced on  $S$ is denoted by $G[S]$, and  $G-S:=G[V(G)\setminus S]$. 
 If $F\subseteq E(G)$, then $G-F$ is obtained from $G$ by deleting all
 the edges of $F$. 
Let $V_1,
V_2\subseteq V(G)$ be two disjoint vertex sets. Then $E_G(V_1,V_2)$ is the set
of edges in $G$  with one end in $V_1$ and the other end in $V_2$, and  $e_G(V_1,V_2):=|E_G(V_1,V_2)|$.  We write $E_G(v,V_2)$ and $e_G(v,V_2)$
if $V_1=\{v\}$ is a singleton.  
Define $\mu(G)=\max\{e_G(u,v)\,:\, u,v\in V(G)\}$
to be the multiplicity of $G$. 
We also write $G[V_1,V_2]$
to denote the bipartite subgraph of $G$ with vertex set $V_1\cup V_2$
and edge set $E_G(V_1, V_2)$. 

For two integers $p,q$, let $[p,q]=\{ i\in \mathbb{Z} \,:\, p \le i \le q\}$. 
 An {\it edge $k$-coloring\/} of a multigraph $G$ is a mapping $\varphi$ from $E(G)$ to the set of integers
$[1,k]$, called {\it colors\/}, such that  no two adjacent edges receive the same color with respect to $\varphi$.  
The {\it chromatic index\/} of $G$, denoted $\chi'(G)$, is defined to be the smallest integer $k$ so that $G$ has an edge $k$-coloring.  
We denote by $\CC^k(G)$ the set of all edge $k$-colorings of $G$.  

In the 1960's, Gupta~\cite{Gupta-67}  and, independently, Vizing~\cite{Vizing-2-classes}  proved
 that for all graphs $G$,  $\Delta(G) \le \chi'(G) \le \Delta(G)+1$. 
This 
leads to a natural classification of simple  graphs. Following Fiorini and Wilson~\cite{fw},  we say a graph $G$ is of {\it class 1} if $\chi'(G) = \Delta(G)$ and of \emph{class 2} if $\chi'(G) = \Delta(G)+1$.  Holyer~\cite{Holyer} showed that it is NP-complete to determine whether an arbitrary graph is of class 1.  
Nevertheless, if $|E(G)|>\Delta(G) \lfloor |V(G)|/2\rfloor$,  then we must use $(\Delta(G)+1)$ colors to edge color $G$. Such graphs are  \emph{overfull}.  An overfull subgraph $H$ of $G$ with  $\Delta(H)=\Delta(G)$
is called a \emph{$\Delta(G)$-overfull subgraph} of $G$. 
 A number of long-standing conjectures listed in {\it Twenty Pretty Edge Coloring Conjectures} in~\cite{StiebSTF-Book} lie in deciding when a 
 graph is overfull.   Chetwynd and  Hilton~\cite{MR848854,MR975994},  in 1986, proposed the following 
conjecture. 
\begin{CON}[Overfull conjecture]\label{overfull-con}
	Let $G$ be a simple graph  with $\Delta(G)>\frac{1}{3}|V(G)|$. Then $\chi'(G)=\Delta(G)$  if and only if $G$ contains no $\Delta(G)$-overfull subgraph.  
\end{CON}

The $3$-critical graph $P^*$, obtained from the Petersen graph by deleting one vertex, has $\chi'(P^*)=4$, 
satisfies $\Delta(P^*)=\frac{1}{3}|V(P^*)|$ but contains no  3-overfull subgraph. 
Thus the degree condition  $\Delta(G)>\frac{1}{3}|V(G)|$ in the conjecture above is best possible.  Applying Edmonds' matching polytope theorem, Seymour~\cite{seymour79}  showed  that whether a graph  $G$ contains an overfull subgraph of maximum degree $\Delta(G)$ can be determined in polynomial time. Thus if the overfull conjecture is true, then the NP-complete problem of 
determining the chromatic index becomes  polynomial-time solvable 
for graphs $G$ with $\Delta(G)>\frac{1}{3}|V(G)|$.
Despite its importance, very little is known about the truth of the overfull conjecture.
It was confirmed only for  graphs with $\Delta(G) \ge |V(G)|-3$ by 
Chetwynd and   Hilton~\cite{MR975994} in 1989. By restricting the minimum degree, 
Plantholt~\cite{MR2082738} in 2004 showed that the overfull conjecture is affirmative  for 
graphs $G$ with  even order $n$ and minimum degree $\delta \ge \sqrt{7}n/3\approx 0.8819 n$. 
The 1-factorization conjecture, which in 2013 
was confirmed for large graphs by Csaba, K\"uhn, Lo, Osthus and Treglown~\cite{MR3545109}, is a special case of the overfull conjecture.
The overfull conjecture was also confirmed for dense quasirandom graphs~\cite{GKO,2104.06253}. 
 In this paper, in supporting the 
 overfull conjecture as well as generalizing the 1-factorization conjecture in an asymptotic way, 
 we obtain the result below.

\begin{THM}\label{thm:1}
For all $0<\ve <1$, there exists $n_0$
such that the following statement holds:
if $G$ is a graph on $2n\ge n_0$ vertices with $\delta(G) \ge (1+\ve)n$,  then $\chi'(G)=\Delta(G)$ if and only if $G$ contains no $\Delta(G)$-overfull subgraph. Furthermore, there is a polynomial time algorithm that finds an optimal coloring. 
\end{THM}

Define $V_i(G)=\{v\in V(G): d_G(v)=i\}$, and we write $V_i$ for $V_i(G)$ if $G$ is clear. 
Furthermore,  $V_{\delta(G)}$ and $V_{\Delta(G)}$ are simply written as  $V_\delta$
and $V_\Delta$, respectively. 
The proof of Theorem~\ref{thm:1} is based on the following result. 

\begin{THM}\label{thm:D-coloring}
	For all $0<\ve <1$, there exists $n_0$
	such that the following statement holds. 
	If $G$ is a graph on $2n\ge n_0$ vertices satisfying  one of the following three conditions:
	\begin{enumerate}[(a)]
		\item $G$ is regular with $\delta(G) \ge (1+4\ve /5)n $,
		\item $G$ has two distinct vertices $x,y$ such that $d(x)=d(y) \ge (1/2+3\ve/2)n$, for all $z\in V(G)\setminus\{x,y\}$,  $d(z)=\Delta(G) \ge (1+\ve)n$, and $\Delta(G)-\delta(G)\le (1/2-\ve/2)n$, 
 		\item $\Delta(G)-\delta(G) \ge n^{6/7}$, $|V_\delta|  \ge n^{6/7}$ and $|V_\Delta|\ge n+1$, and $\delta(G) \ge (1+\ve)n$, 
	\end{enumerate}
	then $\chi'(G)=\Delta(G)$.    
	Furthermore, there is a polynomial time algorithm
	that finds an optimal coloring.     
\end{THM}

The remainder of this paper is organized as follows.
In the next section, we introduce some notation and   preliminary results.
In   Section 3, we prove 
Theorem~\ref{thm:1} by applying Theorem~\ref{thm:D-coloring}. 
Theorem~\ref{thm:D-coloring}
is then proved in the last section.

\section{Notation and preliminaries}

Let $G$ be a multigraph and 
$\varphi\in \CC^k(G)$ for  some integer $k\ge 0$. 
For any $v\in V(G)$, the set of colors \emph{present} at $v$ is 
$\varphi(v)=\{\varphi(e)\,:\, \text{$e$ is incident with $v$} \}$, and the set of colors \emph{missing} at $v$ is $\pbar(v)=[1,k]\setminus\varphi(v)$.
For a subset $X$ of $V(G)$ and a color $i\in [1,k]$, define 
$\pbar_X^{-1}(i)= \{v\in X: i\in \pbar(v)\}$,   
and we  write $\pbar^{-1}(i)$ for $\pbar_{V(G)}^{-1}(i)$. 
An edge $k$-coloring of a multigraph $G$ is said to be \emph{equalized} if each color
class contains either $\lfloor |E(G)|/k \rfloor$ or $\lceil |E(G)|/k \rceil$ edges.  

 For $x\in V(G)$, the \emph{deficiency} of $x$
in $G$ is $\df_{G}(x):=\Delta(G)-d_G(x)$. For $X\subseteq V(G)$, 
$\df_{G}(X)=\sum_{x\in X} \df_{G}(x)$.  
We simply write $\df_{G}(V(G))$
 as $\df(G)$.
A subgraph $H$ of $G$ with an odd order is \emph{$\Delta(G)$-full} if $|E(H)|=\Delta(G)\lfloor |V(H)|/2\rfloor$.

We will use the following notation: $0<a \ll b \le 1$. 
Precisely, if we say a claim is true provided that $0<a \ll b \le 1$, 
then this means that there exists a non-decreasing function $f:(0,1]\rightarrow (0,1]$ such that the statement holds for all $0<a,b\le 1$ satisfying $a \le f(b)$. 

In the 1960's, Gupta~\cite{Gupta-67} and, independently, Vizing~\cite{Vizing-2-classes}   provided an upper bound on the chromatic index
of multigraphs, and K\"onig~\cite{MR1511872} gave an exact value 
of the chromatic index for bipartite multigraphs. 

\begin{THM}[\cite{Gupta-67, Vizing-2-classes}]\label{chromatic-index}
Every multigraph  $G$ satisfies $\chi'(G) \le \Delta(G)+\mu(G)$. 
\end{THM}

\begin{THM}[\cite{MR1511872}]\label{konig}
	Every bipartite multigraph $G$ satisfies $\chi'(G)=\Delta(G)$. 
\end{THM}

%
McDiarmid~\cite{MR300623} observed the following result.
\begin{THM}\label{lem:equa-edge-coloring}
	Let $G$ be a multigraph with chromatic index $\chi'(G)$. Then for all $k\ge \chi'(G)$, there is an equalized edge-coloring of $G$ with $k$ colors. 
\end{THM}

Let $G$ be a multigraph, $k\ge 0$ be an integer and $\varphi \in \CC^k(G)$.
There is a polynomial time algorithm to modify $\varphi$ into
an  equalized edge-coloring of $G$ with $k$ colors. To see this, 
suppose $\varphi$
is not equalized and so we take two colors $i,j\in[1,k]$
such that   $\left||\pbar^{-1}(i)|-|\pbar^{-1}(j)|\right|$ 
is largest.  Since $\varphi$
is not equalized,  $\left||\pbar^{-1}(i)|-|\pbar^{-1}(j)|\right| \ge 4$. 
Assume by symmetry that $|\pbar^{-1}(i)|-|\pbar^{-1}(j)| \ge 4$. 
Consider the submultigraph of $G$ induced on the set of edges colored 
by $i$ or $j$, then the submultigraph must have a component that is a path $P$
starting at an edge colored by  $j$ and ending at an edge colored by $j$. 
By swapping the colors $i$ and $j$ along this path $P$, 
we decreased $|\pbar^{-1}(i)|-|\pbar^{-1}(j)|$ by 4. 
Repeating this process, we can obtain an equalized edge-coloring of $G$ with $k$ colors after at most $k^2 |V(G)|$ rounds. 


Given an edge coloring of $G$ and a given color $i$, 
since  vertices presenting $i$
are saturated by the matching consisting of all edges colored by $i$, we have the Parity Lemma below. The result had appeared in many papers, for example, see~\cite[Lemma 2.1]{MR2028248}.

\begin{LEM}[Parity Lemma]
	Let $G$ be a multigraph and $\varphi\in \CC^k(G)$ for some integer $k\ge \Delta(G)$. 
	Then 
	$|\pbar^{-1}(i)| \equiv |V(G)| \pmod{2}$ for every color $i\in [1,\Delta(G)]$. 
\end{LEM}

We need the following classical result of Hakimi~\cite{MR148049} on multigraphic degree sequence. 
\begin{THM}\label{thm:degree-seq}
	Let $0 \le d_n \le \ldots \le d_1$ be integers. Then there exists a multigraph $G$
	on vertices $x_1,\ldots, x_n$ such that $d_G(x_i)=d_i$
	for all $i$ if and only if $\sum_{i=1}^nd_i$ is even and $\sum_{i>1}d_i \ge d_1$. 
\end{THM}

Though it is not explicitly stated in~\cite{MR148049}, the inductive proof yields a polynomial
time algorithm which finds an appropriate multigraph if it exists.

\begin{THM}[\cite{MR47308}]\label{thm:Dirac}
Let  $G$ be a graph on $n\ge 3$ vertices. If $\delta(G) \ge \frac{n}{2}$, then $G$ is hamiltonian; and if    $\delta(G) \ge \frac{n+1}{2}$, then $G$ is hamiltonian-connected. 
\end{THM}

Following the proof of Dirac~\cite{MR47308}, a hamiltonian cycle can be constructed in polynomial time in $n$
if $\delta(G) \ge \frac{n}{2}$. In fact, there is a polynomial time algorithm that constructs the closure of a graph $G$ and finds a hamiltonian cycle of $G$
if its closure is a complete graph (see~\cite[Exercise 4.2.15, page 62]{MR0411988}).

\begin{LEM}\label{lem:numberofD}
	Let $G$ be an $n$-vertex simple graph such that all  vertices of degree less than $\Delta(G)$ are mutually adjacent in $G$. Then $|V_\Delta|> \frac{n}{2}$. 
\end{LEM}

\pf Suppose the set $A$ of maximum degree vertices has cardinality $k$, and the number of vertices of degree less than maximum degree is $k + r$ with $r \geq 0$.  Deleting $r$ vertices not in A, we get  a new graph  $H$ with $2k$ vertices, $k$ of them forming $A$, and the remaining $k$ forming a set of vertices $B$ such that each vertex in $B$ has degree less than each vertex of $A$ in $H$. But $B$ induces a complete graph in $H$ so in $H$ the sum of the vertex degrees in $A$ is less than or equal to the degree sum of the vertices in $B$.  
Since every vertex of $V(G)\setminus V(H)$ is adjacent in $G$ to every vertex of $B$,
it follows that in $G$ the sum of the vertex degrees in $A$ is less than or equal to the degree sum of the vertices in $B$.
This gives a contradiction. 
\qed 


The  two lemmas below concern  existences of overfull subgraphs in simple graphs. 

\begin{LEM}[\cite{MR1814514}]\label{lem:one-vertex-delted-overfull}
	Let  $G$  be a simple graph of  even order  $n$ with $\delta(G)>\frac{n}{2}$.  If  $H$  is an induced proper subgraph of  $G$  such that  $H$  is either  $\Delta(G)$-overfull or $\Delta(G)$-full, then  $H=G-v$  for some vertex  $v\in V_\delta$.
\end{LEM}

\begin{LEM}\label{lem:overfull2}
	Let  $G$  be a simple graph of  even order  $n$ with $\delta(G)>\frac{n}{2}$.  Then $G$
	contains no $\Delta(G)$-overfull subgraph if $|V_\delta| \ge 2$. 
\end{LEM}

\pf  Let $x,y\in V_\delta$ be distinct. Then 
$\sum_{v\in V(G-x)}(\Delta(G)-d_{G-x}(v))=d_G(x)+(\Delta(G)-d_G(y))+ \df_G(V(G)\setminus \{x,y\}) \ge \Delta(G)$. 
Thus $G-x$ is not $\Delta(G)$-overfull. By Lemma~\ref{lem:one-vertex-delted-overfull}, $G$
contains no  $\Delta(G)$-overfull subgraph. 
\qed 

%
%


\begin{LEM}\label{lem:D-full}
	Let $0<\ve <1$, $n_0$ be a positive integer, 
	and $G$ be a graph on $2n\ge n_0$ vertices with $\delta(G) \ge (1+\ve)n$.  If  $G$ contains a $\Delta(G)$-full subgraph, then $G$ contains a spanning $\delta(G)$-regular subgraph obtained from $G$ by deleting $\Delta(G)-\delta(G)$ matchings iteratively. 
	As a consequence, $\chi'(G)=\Delta(G)$. Furthermore, there is a polynomial time algorithm that finds an optimal coloring. 
\end{LEM}

\pf Define $g=\Delta(G)-\delta(G)$.  If $G$ is regular, then
we are done by Theorem~\ref{thm:D-coloring}. 
Thus $G$
is not regular and so $g\ge 1$.  The graph $G$ contains a $\Delta(G)$-full subgraph, 
which by Lemma~\ref{lem:one-vertex-delted-overfull}  must be $G-x$
for some vertex $x\in V_\delta$. 
Also, if $G$ contains a $\Delta(G)$-overfull subgraph,
then $G-x$ must be $\Delta(G)$-overfull also by Lemma~\ref{lem:one-vertex-delted-overfull}. Since  $G-x$ is $\Delta(G)$-full, we conclude that  
$G$ contains no $\Delta(G)$-overfull subgraph and so has another  vertex of degree less than $\Delta(G)$.   
We let 
$y\in V(G)\setminus\{x\}$ such that $d_G(y)$ is smallest among all vertices in $V(G)\setminus \{x\}$. 
Since $G-x$ is $\Delta(G)$-full, we have  $\Delta(G)=\df(G-x)=d_G(x)+(\Delta(G)-d_G(y))+\df_G(V(G)\setminus \{x,y\})$.
As $d_G(x)=\delta(G)$, if $d_G(y)=\delta(G)$, then $\df_G(V(G)\setminus \{x,y\})=0$. This implies that if $d_G(y)=\delta(G)$,  then every vertex from 
$V(G)\setminus\{x,y\}$ has degree $\Delta(G)$ in $G$; and if $d_G(y)>\delta(G)$,
then as $y$ is chosen to have smallest degree in $G$ among vertices from $V(G)\setminus \{x\}$,  $V(G)\setminus\{x,y\}$ contains no vertex of degree $\delta(G)$ in $G$. Since $\delta(G-x-y) \ge \frac{n}{2}-1$,  $G-x-y$
has a hamiltonian cycle by Theorem~\ref{thm:Dirac}. 
As $n-2$
is even, we know that $G-x-y$ has a perfect matching $M_1$. 
Now we have $\delta(G-M_1)=\delta(G)$ and $\Delta(G-M_1)-\delta(G-M_1)=g-1<g$. 
Let $G_1=G-M_1$. Since 
\begin{eqnarray*}
	\df(G_1-x)&=&d_G(x)+(\Delta(G_1)-d_G(y))+\df_{G_1}(V(G)\setminus \{x,y\})\\
	&=&d_G(x)+\df_G(V(G)\setminus\{x\})-1\\
	&=& \df(G-x)-1=\Delta(G)-1=\Delta(G_1), 
\end{eqnarray*}
we see that 
$G_1-x$ is $\Delta(G_1)$-full.  Thus  we may repeat the procedure, and 
reach a $\delta(G)$-regular graph  $G^*$ after  taking $g$ matchings 
$M_1,\ldots, M_g$. 

Now by Theorem~\ref{thm:D-coloring}, $\chi'(G^*)=\Delta(G^*)=\delta(G)$. Coloring each of the $g$ matchings 
$M_1, \ldots, M_g$ using a different color together with an edge $\delta(G)$-coloring of $G^*$ gives an edge $\Delta(G)$-coloring of $G$. Thus $\chi'(G)=\Delta(G)$. 

It is polynomial-time to find a hamiltonian cycle in graphs $H$ with $\delta(H) \ge \frac{1}{2}|V(H)|$ by the comments immediately after Theorem~\ref{thm:Dirac}. 
Thus all the matchings $M_1, \ldots, M_g$ can be found in polynomial time. 
As an optimal edge coloring can be found in polynomial time for graphs satisfying the conditions in 
Theorem~\ref{thm:D-coloring}, we can find an edge $\delta(G)$-coloring of $G^*$
in polynomial time. Therefore, there is a polynomial time algorithm
that finds an edge $\Delta(G)$-coloring for $G$. 
 \qed 

\begin{LEM}\label{lem:matching-in-bipartite}
	Let $G[X,Y]$ be bipartite graph with $|X|=|Y|=n$. Suppose $\delta(G)=t$ 
for some $t\in [1,n]$, 
	and except at most $t$ vertices all other vertices of $G$ have degree 
	at least $n/2$ in $G$. Then 
	$G$ has a perfect matching. 
\end{LEM}
\pf We show that $G[X, Y]$ satisfies Hall's Condition. If not, 
we let $S\subseteq X$ with smallest cardinality such that 
$|S|>|N_{G}(S)|$. By this choice, $|S|=|N_{G}(S)|+1$ and $|N_{G}(S)|<|Y|$.  
As  $|S|>|N_{G}(S)|$, it follows that 
$|S| \ge \delta(G)+1 \ge t+1$. As $G$
has at most $t$ vertices of degree less than $n/2$,
it then follows that $|S|>n/2$. 
Thus $|X\setminus S|  <n/2$. Since $|N_{G}(S)|<|Y|$, there exists $y\in Y\setminus N_{G}(S)$ such that $N_{G}(y)\subseteq X\setminus S$. As $\delta(G) \ge t$, we have $|X\setminus S| \ge t$. 
As $|Y\setminus N_{G}(S)|=|Y|-|S|+1 =|X|-|S|+1\ge t+1$ and $G$
has at most $t$ vertices of degree less than $n/2$, $Y\setminus N_{G}(S)$
contains a vertex of degree at least $n/2$ in $G$. 
However $|X\setminus S| <n/2$,
we obtain a contradiction. Hence $G$ has a perfect matching.
\qed 

A path $P$ connecting two vertices $u$ and $v$ is called 
a {\it $(u,v)$-path}, and we write $uPv$ or $vPu$ in specifying the two endvertices of 
$P$. Let $uPv$ and $xQy$ be two disjoint paths. If $vx$ is an edge, 
we write $uPvxQy$ as
the concatenation of $P$ and $Q$ through the edge $vx$. If $P$ is a path and $x,y\in V(P)$,
then $xPy$ is the subpath of $P$ with endvertices $x$ and $y$.

\begin{LEM}\label{lem:path-decomposition}
	Let $0<1/n_0 \ll  \ve <1$, and $G$ be 
	graph on $n\ge n_0$ vertices such that $\delta(G) \ge (1+\ve) n/2$. 
	Moreover, let $M=\{a_1b_1,\ldots,a_tb_t\}$ be a matching in the complete graph on $V(G)$ of size at most $\ve n/8$. 
	Then there exist vertex-disjoint path $P_1,\ldots, P_t$ in $G$ such that $\bigcup V(P_i)=V(G)$
	and $P_i$ joins $a_i$ to $b_i$,  and these paths can be found in polynomial time.
\end{LEM}
\pf For $i\in [1,t-1]$, $|N_G(a_i)\cap N_G(b_i)| \ge \ve n$, 
so we can greedily find  vertices $c_i\in N_G(a_i)\cap N_G(b_i)$
such that $c_i\ne c_j$ for distinct $i,j\in [1,t-1]$. Thus we let $P_i=a_ic_ib_i$. 
Let $G^*=G-\bigcup_{i=1}^{t-1} V(P_i)$. Then $\delta(G^*)  \ge (1+\ve )n/2-3(t-1) \ge (1+\ve/8 )n/2$, and so $G^*$ is hamiltonian-connected by Theorem~\ref{thm:Dirac}. Thus  we can find an $(a_t,b_t)$-hamiltonian path  $P_t$ in $G^*$. 

It is clear that each of $P_1,\ldots, P_{t-1}$ can be found in polynomial time. 
For the path $P_t$, we construct it as below. By the comments immediately after Theorem~\ref{thm:Dirac}, we  can find a hamiltonian cycle  $C$ of 
$G^*$ in polynomial time.  By taking a longer segment  between $a_t$ and $b_t$ from $C$, we get in $G^*$
an $(a_t,b_t)$-path  $Q_1$ that contains at least $|V(G^*)|/2$ vertices. We will 
extend $Q_1$ into a hamiltonian $(a_t,b_t)$-path of $G^*$. 
Denote by $Q_2$ the remaining segment of $C$ that is disjoint from $Q_1$
and let $c$ and $d$ be the endvertices of $Q_2$. Let $|V(Q_2)|=p$.
Then as $\delta(G^*)\ge (1+\ve/8 )n/2$, each of $c$ and $d$
has on $Q_1$ at least $(1+\ve/8 )n/2-(p-1)=(1+\ve/8 )n/2-p+1$ neighbors. Since $2((1+\ve/8 )n/2-p)+p+|V(Q_2)|>|V(G^*)|$, 
it follows that one of the following two situations must happen:
(a) there is  a vertex $c_1 \in N_{G^*}(c)\cap V(Q_1)$ and a vertex $d_1 \in N_{G^*}(d)\cap V(Q_1)$ such that $c_1Q_1d_1$ contains less than $p+2$ vertices, and (b)  $c$ or $d$ has on $Q_1$ two neighbors that are consecutive on $Q_1$. When (a) happens,  
assume by symmetry that $c_1$ is between $a_t$ and $d_1$ on $Q_1$, 
then  $Q_1^*=a_tQ_1c_1cQ_2dd_1Q_1b_t$ is longer than $Q_1$ and the component of $G^*-V(Q_1^*)$ still contains a hamiltonian path.  Similarly, we can extend $Q_1$
into a longer $(a_t,b_t)$-path such that the subgraph of $G^*$
outside the path is hamiltonian if (b) happens. Repeating this procedure at most $n/2$ times, we obtain a hamiltonian $(a_t,b_t)$-path of $G^*$. 
Therefore, all the path $P_1, \ldots, P_t$ can be found in polynomial time.  
\qed

\section{Proof of Theorem~\ref{thm:1}}


\begin{THM1}
	For all $0<\ve <1$, there exists $n_0$
	such that the following statement holds:
	if $G$ is a graph on $2n\ge n_0$ vertices with $\delta(G) \ge (1+\ve)n$,  then $\chi'(G)=\Delta(G)$ if and only if $G$ contains no $\Delta(G)$-overfull subgraph. Furthermore, there is a polynomial time algorithm that finds an optimal coloring. 
\end{THM1}

\pf 
Choose constants $\ve$ and positive integer $n_0$  such that $0<1/n_0 \ll  \ve$.

If $G$ is regular, then we are done by Theorem~\ref{thm:D-coloring}. 
Thus we assume that $G$ is not regular. If $G$ contains a $\Delta(G)$-overfull subgraph,
then $\chi'(G)=\Delta(G)+1$. Thus we assume that $G$ contains no $\Delta(G)$-overfull subgraph. As a consequence, $\df(G) \ge \Delta(G)$. 
By Lemma~\ref{lem:D-full}, we may assume that  $G$
contains no $\Delta(G)$-full subgraph.   
 Therefore, if two vertices with degree less than  $\Delta(G)$  are not adjacent in $G$, we may add the edge between them without creating a  $\Delta(G)$-overfull subgraph, or increasing $\Delta(G)$. We iterate this edge-addition procedure.  
 If  at some point we create a $\Delta(G)$-full subgraph, the result follows by Lemma~\ref{lem:D-full}.  Otherwise, we reach a point where we may now assume that in  $G$  all vertices with degree less than $\Delta(G)$ are mutually adjacent, and  so by Lemma~\ref{lem:numberofD}, we have $|V_\Delta| \ge n+1$.

Define $n_1=|V_\delta|$. Note that $n_1<n$. 
If $n_1 \ge  n^{6/7}$ and $\Delta(G)-\delta(G) \ge n^{6/7}$, then   we are done by Theorem~\ref{thm:D-coloring}. Thus we assume $n_1 < n^{6/7}$ or $\Delta(G)-\delta(G) < n^{6/7}$, and we consider the two cases below. 
We call a  vertex of degree less than $\Delta(G)$ but greater than $\delta(G)$ a \emph{middle degree}  vertex. 

{\bf \noindent Case 1. $n_1< n^{6/7} $. }

Note that for any $v\in V(G)\setminus V_\delta$,  $\delta(G-v-V_\delta) \ge n$
and so both $G-V_\delta$ and $G-v-V_\delta$ are hamiltonian by Theorem~\ref{thm:Dirac}. 
Thus if $G-V_\delta$ and $G-v-V_\delta$ have even order, then they each have a perfect matching. 
Hence if $n_1$ is even, we can decrease $\Delta(G)-\delta(G)$ but preserve $\delta(G)$
in deleting a perfect matching  $M$ of $G-V_\delta$. If $n_1$ 
is odd but $G$ has a middle degree vertex $v$, we can decrease $\Delta(G)-\delta(G)$ but preserve $\delta(G)$
in deleting a perfect matching  $M$ of $G-v-V_\delta$.  Denote by $G_1$
the reduced graph from $G$ by deleting $M$ in either of these two cases. 
If $|V_\delta| \ge 2$, then  as $V_\delta\subseteq V_\delta(G_1)$, we know that $G_1$
still contains no $\Delta(G_1)$-overfull subgraph by Lemma~\ref{lem:overfull2}. 
Thus $|V_\delta|=1$. Let  $V_\delta=\{u\}$.  Note that $u \in V_\delta(G_1)$. 
Then $\df(G_1-u) =d_G(u)+(\Delta(G_1)-d_G(v))+\df_{G_1}(V(G_1)\setminus \{u,v\})=d_G(u)+(\Delta(G)-1-d_G(v))+\df_{G}(V(G)\setminus \{u,v\})=\sum_{w\in V(G-u)}(\Delta(G)-d_{G-u}(w))-1$. Since $G$
contains no $\Delta(G)$-overfull subgraph, we have 
$ \sum_{w\in V(G-u)}d_{G-u}(w) \ge \Delta(G)$.  Thus $\df(G_1-u) \ge \Delta(G)-1=\Delta(G_1)$
and so $G_1$ contains no $\Delta(G_1)$-overfull subgraph by Lemma~\ref{lem:one-vertex-delted-overfull}. 
Furthermore, $\chi'(G_1)=\Delta(G_1)$ implies that $\chi'(G)=\Delta(G)$. 
Thus in these two cases, we can consider $G_1$ in place of $G$ and show that $G_1$
is a class 1 graph.

Thus we assume  $n_1$ is odd and $G$ has no middle degree vertex. This in particular, implies that $\delta(G)$ and $\Delta(G)$ have the same parity. 
As $G$ has no $\Delta(G)$-overfull subgraph, $|V_\delta|\ge 3$. 
Let $x,y\in V_\delta$ be distinct.  We find a perfect matching $M_{11}$
in $G-(V_\delta\setminus \{x\})$ and a perfect matching $M_{12}$
in $G-(V_\delta\setminus \{y\})$. The matchings exist by Theorem~\ref{thm:Dirac}.    Let $G_1=G-M_{11}-M_{12}$. 
We repeat this same process and find a perfect matching $M_{21}$
in $G_1-(V_\delta\setminus \{x\})$ and a perfect matching $M_{22}$
in $G_1-(V_\delta\setminus \{y\})$. For $i\in [2,(\Delta(G)-\delta(G))/2]$, we let $G_i=G_{i-1}-M_{i1}-M_{i2}$. 
We have $d_{G_i}(x)=d_{G_i}(y)=\delta(G)-i$. 
As $\Delta(G)-\delta(G) \le 2n-(1+\ve)n=(1-\ve)n$, we see that 
$d_{G_i}(x)=d_{G_i}(y)=\delta(G)-i \ge (1+\ve)n-\frac{1}{2}(1-\ve)n\ge (1/2+3\ve/2)n$. 
For any vertex $z\in V(G_i)\setminus \{x,y\}$, $d_{G_i}(z) \ge \delta(G) \ge (1+\ve)n$. 
Let $x^*$ be a neighbor of $x$ in $G_i-(V_\delta\setminus \{x\})$ 
and $y^*$ be a neighbor of $y$ in $G_i-(V_\delta\setminus \{y\})$. 
Then  $G_i-(V_\delta\setminus \{x\})-\{x,x^*\}$ has a perfect matching $M_{(i+1)1}^*$, 
and $G_i-(V_\delta\setminus \{y\})-\{y,y^*\}$ has a perfect matching $M_{(i+1)2}^*$. 
Let $M_{(i+1)1}=M_{(i+1)1}^* \cup \{xx^*\}$ and $M_{(i+1)2}=M_{(i+1)2}^* \cup \{yy^*\}$. 
Thus for each $i\in [2,(\Delta(G)-\delta(G))/2]$, we find matchings $M_{i1}$
and $M_{i2}$ respectively 
from $G_{i-1}-(V_\delta\setminus \{x\})$ and $G_{i-1}-(V_\delta\setminus \{y\})$. 

We claim $G^*:=G_{(\Delta(G)-\delta(G))/2}$ satisfies Condition (b) 
of Theorem~\ref{thm:D-coloring}.  By the analysis above, we have $d_{G^*}(x)=d_{G^*}(y) \ge (1/2+3\ve/2)n$, $d_{G^*}(z)=\Delta(G^*) =\delta(G)\ge (1+\ve)n$ for all $z\in V(G^*)\setminus\{x,y\}$.  
Also $\Delta(G^*)-\delta(G^*)=\delta(G)- (\delta(G)-(\Delta(G)-\delta(G))/2)=\frac{1}{2}(\Delta(G)-\delta(G))
\le \frac{1}{2}( 2n-(1+\ve)n)=(1/2-\ve/2)n$.  
By Theorem~\ref{thm:D-coloring}, $\chi'(G^*)=\Delta(G^*)=\delta(G)$. 
Taking an edge $\delta(G)$-coloring of $G^*$,  coloring edges  in $M_{i1}$  with color $\delta(G)+2i-1$
and coloring edges in $M_{i2}$  with color $\delta(G)+2i$ for each $i\in [1, (\Delta(G)-\delta(G))/2]$, 
we obtain an edge $\Delta(G)$-coloring of $G$.

{\bf \noindent Case 2. $\Delta(G)-\delta(G) <n^{6/7}$. }

Let  $V(G)=\{x_1,\ldots,x_{2n}\}$ and we assume $\df_G(x_1) \ge \ldots \ge  \df_G(x_{2n})=0$. 
Since $x_1$ has the smallest degree in $G$ and $G-x_1$ is not  $\Delta(G)$-overfull by our assumption, 
$\sum_{i\ge 2}\df_{G}(x_i) \ge \df_{G}(x_1)$. Since $|V(G)|=2n$
is even,  $\sum_{i\ge 1}\df_{G}(x_i) $ is even. 
Then by Theorem~\ref{thm:degree-seq}, there exists a multigraph $H$
on $V(G)$ such that $d_H(x)=\df_{G}(x_i)$
 for each $i\in [1,2n]$. This multigraph $H$
will aid us to find a spanning regular subgraph of $G$. 

Note that $\Delta(H)=\df_{G}(x_1)=\Delta(G)-\delta(G) < n^{6/7}$ and $H$ contains isolated vertices.
Thus $\chi'(H) \le \Delta(H)+\mu(H) \le 2\Delta(H) \le 2n^{6/7}$. Hence we 
can greedily partition $E(H)$ into $k\le 5n^{6/7}/\ve$ matchings $M_1,\ldots, M_k$ 
each of size at most $\ve n/5$.  Now we take out linear forests  
from $G$ by applying Lemma~\ref{lem:path-decomposition} with $M_1,\ldots, M_k$. 
More precisely, define spanning subgraphs $G_0, \ldots, G_k$
of $G$ and edge-disjoint linear forests $F_1, \ldots, F_k$
such that 
\begin{enumerate}[(1)]
	\item $G_0:=G$ and $G_i=G_{i-1}-E(F_i)$ for $i\in [1,k]$,
	\item $F_i$ is a spanning linear forest (each vertex of $G_{i-1}$ has degree 1 or 2 in $F_i$) in $G_{i-1}$ whose leaves are precisely the 
	vertices in $M_i$. 
\end{enumerate}

Let $G_0=G$ and suppose that for some $i\in[1,k]$, we already defined $G_0, \ldots, G_{i-1}$
and $F_1, \ldots, F_{i-1}$. As $\Delta(F_1\cup \ldots \cup F_{i-1}) \le 2(i-1) \le 10 n^{6/7}/\ve$,
it follows  that $\delta(G_{i-1}) \ge (1+\ve) n-10 n^{6/7}/\ve\ge (1+4\ve /5)n$. Since $M_i$ has size at most $\ve n/5$, we can apply Lemma~\ref{lem:path-decomposition} to $G_{i-1}$ and $M_i$ and obtain a spanning linear forest  $F_{i}$ in $G_{i-1}$ whose leaves are precisely the vertices in $M_i$. 
Set $G_i:=G_{i-1}-E(F_{i})$. 

We claim that $G_k$ is regular. Consider any vertex $u\in V(G_k)$. 
For every $i\in[1,k]$, $d_{F_i}(u)=1$ if $u$ is an endvertex of some edge of $M_i$
and $d_{F_i}(u)=2$ otherwise. 
Since $M_1, \ldots, M_k$ partition $E(H)$, we know that $\sum\limits_{i=1}^{k}d_{F_i}(u)=2k-d_H(u)=2k-\df_{G}(u)$.
Thus 
$$d_{G_k}(u)=d_{G}(u)-\sum\limits_{i=1}^{k}d_{F_i}(u)=d_{G}(u)-(2k-\df_{G}(u))=\Delta(G)-2k.$$
Note that $\Delta(G) \ge (1+\ve)n-10 n^{6/7}/\ve\ge (1+4\ve /5)n$. 
Now $\chi'(G_k)=\Delta(G_k)$ by Theorem~\ref{thm:D-coloring}. 
We color  the edges of $F_i$ 
using 2 distinct colors from $[\Delta(G)-2k+1, \Delta(G)]$  for each $i\in[1,k]$. 
It is clear that any edge $\Delta(G_k)$-coloring of $G_k$ together with 
this coloring of  $\bigcup_{i=1}^kF_i$ gives an edge coloring 
of $G$ using $\Delta(G_k)+2k=\Delta(G)$ colors.

We lastly check that the procedure above yields a polynomial time algorithm. Given $G$, 
taking a vertex $u$ of minimum degree in $G$,  we first check if 
$G-u$ is $\Delta(G)$-overfull. If yes, then $\chi'(G)=\Delta(G)+1$ and 
$G$ can be edge colored using $\Delta(G)+1$ colors in polynomial
time~\cite{MR1156837}.   Thus $G$ contains no $\Delta(G)$-overfull 
subgraph. If $G$ contains a $\Delta(G)$-full subgraph, then an 
edge $\Delta(G)$-coloring of $G$ can be found in polynomial time by 
Lemma~\ref{lem:D-full}. Thus $G$ contains no $\Delta(G)$-full subgraph. 
If there exist nonadjacent $u,v\in V(G)\setminus V_\Delta$, we add the edge 
$uv$ in $G$. If we reach a point where the resulting graph contains 
a $\Delta(G)$-full subgraph, we then find an 
edge $\Delta(G)$-coloring  of the graph  in polynomial time by 
Lemma~\ref{lem:D-full}, which also gives an edge $\Delta(G)$-coloring of 
$G$. Thus we assume that every two  vertices from $V(G)\setminus V_\Delta$ 
are adjacent in $G$. If $G$ is in Condition (c) of Theorem~\ref{thm:D-coloring},
then we find an 
edge $\Delta(G)$-coloring of $G$  in polynomial time by Theorem~\ref{thm:D-coloring}. 
Thus we have Case 1 or Case 2 as described in this proof. 
If $G$ is in Case 1, it is polynomial time to find the desired matchings (basically find hamiltonian cycles of even length in graphs with large minimum degree by  the comments immediately after Theorem~\ref{thm:Dirac})
to reduce $G$ into a graph satisfying one of the conditions in Theorem~\ref{thm:D-coloring}. Then  we find an 
edge $\Delta(G)$-coloring of $G$  in polynomial time by Theorem~\ref{thm:D-coloring}. 
If $G$ is in Case 2, 
then 
can construct an edge $\Delta(G)$-coloring of $G$ through the process as described in Case 2. Since Theorem~\ref{thm:degree-seq},  Lemma~\ref{lem:path-decomposition} and Theorem~\ref{thm:D-coloring}  give appropriate
running time statements, this can be achieved in time polynomial in $n$. 
\qed

\section{Proof of Theorem~\ref{thm:D-coloring}}

The  proofs in this section  follow and extend ideas of Vaughan from~\cite{MR2993074}, where the techniques were for regular multigraphs
but we modify them for graphs that are not necessarily regular.  
We will need the following result, which was proved using 
Chernoff bound. 


\begin{LEM}[\cite{2104.06253}, Lemma 3.2]\label{lem:partition}
	There exists a positive integer $n_0$ such that for all $n\ge n_0$ the
	following holds. Let $G$ be a graph on $2n$ vertices, and $N=\{x_1,y_1,\ldots, x_t,y_t\}\subseteq V(G)$. 
	Then $V(G)$ can be partitioned into two  parts 
	$A$ and $B$ satisfying the properties below:
	\begin{enumerate}[(i)]
		\item  $|A|=|B|$;
		\item $|A\cap \{x_i,y_i\}|=1$ for each $i\in [1,t]$;
		\item $| d_A(v)-d_B(v)| \le n^{2/3}-1$ for each $v\in V(G)$, where $d_S(v)=|N_G(v)\cap S|$ for any $S\subseteq V(G)$. 
	\end{enumerate}
Furthermore, one such partition can be constructed in $O(2n^3 \log_2 (2n^3))$-time. 
\end{LEM}

\begin{THM2}
For all $0<\ve <1$, there exists $n_0$
such that the following statement holds. 
If $G$ is a graph on $2n\ge n_0$ vertices satisfying  one of the following three conditions:
\begin{enumerate}[(a)]
	\item $G$ is regular with $\delta(G) \ge (1+4\ve /5)n $,
	\item $G$ has two distinct vertices $x,y$ such that $d(x)=d(y) \ge (1/2+3\ve/2)n$, for all $z\in V(G)\setminus\{x,y\}$,  $d(z)=\Delta(G) \ge (1+\ve)n$, and $\Delta(G)-\delta(G)\le (1/2-\ve/2)n$, 
	\item $\Delta(G)-\delta(G) \ge n^{6/7}$, $|V_\delta|  \ge n^{6/7}$ and $|V_\Delta|\ge n+1$, and $\delta(G) \ge (1+\ve)n$, 
\end{enumerate}
then $\chi'(G)=\Delta(G)$.    
Furthermore, there is a polynomial time algorithm
that finds an optimal coloring.    
\end{THM2}

\pf  
If $G$ is in Condition (a), we let $N=\emptyset$. 
If $G$ is in Condition (b), we let $N=\{x_1,y_1\}$ where $x_1=x$
and $y_1=y$. If $G$
is in Condition (c), we 
 take $2\lfloor(n-|V_\Delta|)/2 \rfloor$ vertices from $V(G)\setminus V_\Delta$ 
and name them as $x_1,y_1, \ldots, x_t, y_t$, where $t:=\lfloor (n-|V_\Delta|)/2 \rfloor$ and  we assume that the first $\lfloor |V_\delta|/2\rfloor$ pairs of vertices $x_i,y_i$ are all from $V_\delta$.
Let $N=\{x_1,y_1, \ldots, x_t, y_t\}$.   
Applying Lemma~\ref{lem:partition} on $G$ and $N$,
we obtain a partition $\{A, B\}$ of $V(G)$ satisfying the following properties:
\begin{enumerate}[P.1]
	\item $|A|=|B|$;
	\item \text{$|A\cap \{x_i,y_i\}|=1$ for each $i\in [1,t]$};
	\item \text{$| d_{A}(v)-d_{B}(v)| \le n^{2/3}-1$ for each $v\in V(G)$}. 
\end{enumerate}
  Thus when $G$ is in Condition (b),  we may assume $x\in A$ and $y\in B$. 
  When $G$ is in Condition (c), we know that $|A\cap V_\delta| \ge \frac{1}{2} (|V_\delta|-1)$, $|B\cap V_\delta| \ge \frac{1}{2} (|V_\delta|-1)$, $|A\cap V_\Delta| \ge n/2$ and $|B\cap V_\Delta| \ge n/2$. 
  By P.3, for any $v\in V(G)$, 
  we have 
  $$
  \frac{1}{2}\left(d_G(v)-n^{2/3}\right) \le d_A(v), d_B(v)  \le  \frac{1}{2}\left(d_G(v)+n^{2/3}\right).
  $$

Let 
$$ G_A=G[A], \quad G_B=G[B], \quad \text{and } \quad H=G[A,B]. $$

To prove the theorem, we will construct an edge coloring 
of $G$ using $\Delta(G)$ colors. 
We provide below an overview of the steps. 
At the start of the process, $E(G)$ is assumed to be uncolored,
and throughout the process, the partial edge coloring of $G$
is always denoted by $\varphi$, which is updating step by step. 

\begin{enumerate}[Step 1]
	\item Define $S=\{v\in V(G)\,:\, \Delta(G)-d_G(v) \ge 7n^{2/3}\}$. 
	Let $k=\max\{\Delta(G_A), \Delta(G_B)\}+1$. By Theorem~\ref{chromatic-index}, we find an edge $k$-coloring  $\varphi$ of  $G_A \cup G_B$. 
	If there exist distinct $u,v\in S\cap A$  or distinct $u,v \in S\cap B$ such that $\pbar(u)\cap \pbar(v)\ne \emptyset$, we add an edge joining $u$
	and $v$ and color the new edge by a color in $\pbar(u)\cap \pbar(v)$. The edge coloring $\varphi$ is updated and we still call it $\varphi$. 
	We iterate this process of adding and coloring edges and call the  multigraphs  resulting from $G_A$ and $G_B$, respectively, $G^*_A$ and $G^*_B$, and call 
	$G^*$ the union of $G^*_A$, $G^*_B$ and $H$. We will modify the current 
	edge coloring, which is still named $\varphi$, such that the following properties are satisfied:
	\begin{enumerate}
		\item[S1.1 ] When $G$ is in Conditions (a) or (b), $$
		|\pbar^{-1}_A(i)|=|\pbar^{-1}_B(i)| \quad \text{for every $i\in [1,k]$}. 
		$$
		When $G$ is in Condition (c), assume by symmetry that $e(G^*_A) \le  e(G^*_B)$, then 
		$$
		|\pbar^{-1}_A(i)| \ge |\pbar^{-1}_B(i)| \quad \text{for every $i\in [1,k]$}. 
			$$ 	\label{s1.1}
		\item[S1.2 ] 
		\begin{equation}\tag{S1.I}\label{eq-vertices-missing-i}
		\begin{split}
		\sum\limits_{u\in A}|\pbar(u)| \le \frac{n}{2}(n^{2/3}+1)+\frac{n}{2}(6n^{2/3})+k \le 4n^{5/3}-2n,\\
		\sum\limits_{u\in B}|\pbar(u)| \le \frac{n}{2}(n^{2/3}+1)+\frac{n}{2}(6n^{2/3})+k \le 4n^{5/3}-2n,\\
		|\pbar^{-1}_A(i)| \le 4 n^{2/3} \quad \text{and} \quad  |\pbar^{-1}_B(i)| \le 4 n^{2/3} \quad  \text{for every  $i\in [1,k]$.} 
		\end{split}
		\end{equation}
		\label{s1.2}
	\end{enumerate}

	\item Modify the partial edge-coloring of $G^*$ obtained in Step 1 by exchanging alternating paths. When this step is finished, each of the $k$ color class will be a 1-factor of $G^*$. During the process of this step,  a few edges of $H$ will be colored and  a few edges of $G^*_A$ and $G^*_B$ will be uncolored. Denote by  $R_A$ and $R_B$, respectively,  the submultigraphs of $G^*_A$
	and $G^*_B$ consisting of the uncolored edges. The two multigraphs $R_A$ and $R_B$  will initially be empty, but 
	one, two or three edges will be added to at least one of them 
	when each time we exchange an alternating path. 
	The conditions below will be satisfied at the completion of this step:
	\begin{enumerate}
		\item [S2.1] The number of uncolored edges in each of $G_A^*$
		and $G_B^*$ is less than $12 n^{5/3}$. 
		When $G$ is in Conditions (a) or (b), $G_A^*$
		and $G_B^*$ have the same number of uncolored edges; and when $G$
		is in Condition (c), the number of uncolored edges in $G_B^*$
		is greater than or equal to the number of uncolored edges in $G_A^*$ (this follows from our assumption that  $e(G_A^*) \le e(G_B^*)$). 
		\label{s21}
		\item [S2.2] $\Delta(R_A)$ and $\Delta(R_B)$ are less than $n^{5/6} + 1.$
		\label{s22}
		\item [S2.3] Define 
		\begin{eqnarray*}
		S_A&=&\{u\in S\cap A\,:\, d_{G_A^*}(u) \le k-2n^{2/3} \},\\
		S_B&=&\{u\in S\cap B\,:\, d_{G_B^*}(u) \le k-2n^{2/3} \}.
		\end{eqnarray*}
	We require
	\begin{enumerate}[S2.3.1]
		\item Every vertex in $V(G^*)\setminus (S_A\cup S_B)$ is incident in $G^*$ with fewer than $2n^{5/6}$ colored edges of $H$.
		\item When $G$ is in Condition (b),  each of the  vertex from $S_A\cup S_B$
		is incident in $G^*$ with fewer than $(\frac{1}{4}-\frac{1}{5}\ve )n$ colored edges of $H$.
		\item  When $G$ is in Condition (c), each of the  vertex from $S_A\cup S_B$
		is incident in $G^*$ with fewer than $(\frac{1}{2}-\frac{1}{3}\ve )n$ colored edges of $H$. 
	\end{enumerate}
		\label{s23}
	\end{enumerate}
	
	\item We will edge color  $R_A$
	and $R_B$ and a few uncolored edges of $H$ using another 
	$\ell$ colors, where $\ell =\lceil 2n^{5/6}  \rceil$.  
	The goal is to ensure that 
	each of these  $\ell$ new color classes obtained at the completion of Step 3 presents at all vertices from $V(G^*)\setminus V_\delta$ while preserving the $k$ 1-factors already obtained through Steps 1 and 2.  
	\item At the start of Step 4, all of the uncolored edges of $G^*$ belong to $H$. Denote by $R$ the subgraph of $G^*$ consisting of the uncolored edges. 
	It will be shown that $\Delta(R)=\Delta(G^*)-k-\ell$. 
	 This subgraph is bipartite, so we can color its edges using $\Delta(G^*)-k-\ell$
	colors by Theorem~\ref{konig}.  
\end{enumerate}

When Step 4 is completed, we obtain an edge coloring of $G^*$ using exactly $\Delta(G^*)$ colors. We now give the details of each step,
and for concepts that were already defined in the outline above,
we will use them directly.

\bigskip 

\bigskip 

\begin{center}
 	Step 1: Coloring $G_A$ and $G_B$ 
\end{center}

Recall $S=\{v\in V(G)\,:\, \Delta(G)-d_G(v) \ge 7n^{2/3}\}$.
 Note that when $G$ is in Condition (a), $S=\emptyset$; when $G$ is in Condition (b), then $S\subseteq\{x,y\}$;
and   when $G$ is in Condition (c), then $V_\delta\subseteq S$.
Following the operations described in the outline of Step 1, 
for the current edge coloring $\varphi$ of $G_A^* \cup G_B^*$,  
the following statement holds:  $\pbar(u)\cap \pbar(v)=\emptyset$
for any two distinct $u,v\in S\cap A$ or any two distinct $u,v\in S\cap B$.  
Therefore, 
 \begin{equation}\tag{S1.II}\label{eq1}
 \sum_{u\in A\cap S}|\pbar(u)|=\sum_{u\in A\cap S} (k-d_{G^*_A}(u)) \le k, \quad   \sum_{u\in B\cap S}|\pbar(u)|=\sum_{u\in B\cap S} (k-d_{G^*_B}(u)) \le k. 
 \end{equation}

We will in the rest of the proof show that $\chi'(G^*)=\Delta(G^*)$,
this is because $G$ is a subgraph of $G^*$ and  $\Delta(G^*)=\Delta(G)$.
The latter is seen as below:    
for any $u\in S\cap A$, we have 
$$
d_{G^*}(u) \le 
k +e_G(u, B)  \le  \frac{1}{2}(\Delta(G)+n^{2/3})+1+\frac{1}{2}(\Delta(G)-7n^{2/3}+n^{2/3}) \le \Delta(G).
$$
Similarly, we have $d_{G^*}(u) \le \Delta(G)$ for any $u\in S\cap B$.  
In particular, if $u\in V_\delta$, as $\Delta(G)-\delta(G) \ge n^{6/7}$, 
we have 
\begin{equation}\tag{S1.III}\label{eqn2}
d_{G^*}(u) \le  \frac{1}{2}(\Delta(G)+n^{2/3})+1+\frac{1}{2}(\Delta(G)-n^{6/7}+n^{2/3})\le \Delta(G)-\frac{1}{3}n^{6/7}. 
\end{equation}

Let $\varphi_{A}$ and $\varphi_B$ be the restrictions of $\varphi$
on $G^*_A$ and $G^*_B$, respectively. 
By Lemma~\ref{lem:equa-edge-coloring} and the comments immediately below the lemma, we modify 
$\varphi_A$ and $\varphi_B$ into  equitable edge $k$-colorings of $G^*_A$ and $G^*_B$, respectively, and still call $\varphi$ the edge $k$-coloring of $G^*_A\cup G^*_B$
consisting of the modifications of $\varphi_A$ and $\varphi_B$. 
Note that under the new colorings, it is possible that $\pbar(u)\cap \pbar(v) \ne \emptyset$
for some distinct  $ u,v\in S\cap A$ or distinct $u,v\in S\cap B $. 
However  the inequalities in~\eqref{eq1} still hold.

When $G$ is in Conditions (a) or (b),  we have $|S| \le 2$  and $e(G_A)=e(G_B)$ by the partition $\{A,B\}$ of $V(G)$.  Since  $|S\cap A|=|S\cap B| \le 1$, it follows that $G^*_A=G_A$ and $G^*_B=G_B$.  Thus $e(G^*_A)=e(G^*_B)$. Since $\varphi_A$ and $\varphi_B$ 
are equitable edge $k$-colorings of $G^*_A$ and $G^*_B$, 
by renaming some color names in $G^*_A$ if necessary, we  assume 
$$
|\pbar^{-1}_A(i)|=|\pbar^{-1}_B(i)| \quad \text{for every $i\in [1,k]$}. 
$$
When $G$ is in Condition (c), by symmetry, we assume $e(G^*_A) \le e(G^*_B)$.
For the same reasoning as above, we assume  
$$
|\pbar^{-1}_A(i)| \ge |\pbar^{-1}_B(i)| \quad \text{for every $i\in [1,k]$}. 
$$
By the Parity Lemma, $|\pbar^{-1}_A(i)|- |\pbar^{-1}_B(i)|$ is even for every $i\in[1,k]$.
Therefore, we have the statement S1.1 as stated in the outline of Step 1.

Next, we verify that every color $i\in [1,k]$ is missing at a small number of vertices.  Property P.2 of the partition $\{A,B\}$ 
implies  $|A\cap V_\Delta| \ge n/2$ and $|B\cap V_\Delta| \ge n/2$,
and each vertex $u\in V_\Delta$ satisfies $|\pbar(u)| \le n^{2/3}+1$, call this Fact 1.  
 By the definition of $S$, for every $u\in V(G^*)\setminus S$, $d_{G^*}(u)=d_G(u)  >\Delta(G^*)-7n^{2/3}$, and Property P.3
of the  partition $\{A,B\}$ implies 
$d_{G_A^*}(u) \ge \frac{1}{2}d_G(u)-n^{2/3}$ for every $u\in A\setminus S$ and $d_{G_B^*}(u) \ge \frac{1}{2}d_G(u)-n^{2/3}$ for every $u\in B\setminus S$. 
Thus $|\pbar(u)| \le k- (\frac{1}{2}d_G(u)-n^{2/3}) <6n^{2/3}$ for every $u\in V(G^*)\setminus S$, call this Fact 2. 
These two facts 
together with the fact in~\eqref{eq1}, give 
$$
\sum\limits_{u\in A}|\pbar(u)| \le \frac{n}{2}(n^{2/3}+1)+\frac{n}{2}(6n^{2/3})+k \le 4n^{5/3}-2n. 
$$
Similarly, 
$$
\sum\limits_{u\in B}|\pbar(u)| \le \frac{n}{2}(n^{2/3}+1)+\frac{n}{2} (6n^{2/3})+k \le 4n^{5/3}-2n. 
$$
Since $\varphi_A$ and $\varphi_B$ 
are equitable edge $k$-colorings of $G^*_A$ and $G^*_B$, we get 
$$
|\pbar^{-1}_A(i)| \le 4 n^{2/3} \quad \text{and} \quad  |\pbar^{-1}_B(i)| \le 4 n^{2/3}.
$$
Therefore, we have the statement S1.2 as stated in the outline of Step 1.

\begin{center}
	Step 2: Extending  existing color classes into  1-factors  
\end{center}

Each of the $k$ color classes obtained in Step 1 will be extended into 
$k$ 1-factors of $G^*$
through exchanging of alternating paths, which  consist of colored edges 
and uncolored edges. The colored edges and uncolored edges of these alternating paths 
are from $G^*_A\cup G^*_B$ and $H$, respectively. Thus 
during the procedure of Step 2, we will uncolor some of the edges of 
$G^*_A$ and $G^*_B$, and will color some of the edges of $H$. 
Recall that $R_A$ and $R_B$ are the submultigraphs of $G^*_A$
and $G^*_B$ consisting of the uncolored edges, which are 
empty initially. 

To ensure Condition S2.2 is satisfied, we say that an edge  $e=uv\in E(G^*_A\cup G_B^*)$ is \emph{good} if $e\not \in E(R_A \cup R_B)$  and  the degree of $u$ and $v$ 
in both $R_A$ and $R_B$ is less than $n^{5/6}$ (actually, note that when  $uv\in E(G_A^*)$, then the degree of $u$
and $v$ is zero in $R_B$  and vice versa). 
Thus  a good edge can be added to $R_A$ or $R_B$ without violating S2.2.

 By S1.\ref{s1.1},  for each color 
 $i\in [1,k]$, we pair up each vertex from $\pbar^{-1}_B(i)$ 
 with a vertex from $\pbar^{-1}_A(i)$,  and then pair up 
 the remaining unpaired vertices from  $\pbar^{-1}_A(i)$ 
 as $|\pbar^{-1}_A(i)|-|\pbar^{-1}_B(i)|$  is even and we assumed 
 $|\pbar^{-1}_A(i)| \ge |\pbar^{-1}_B(i)|$. Each of those pairs is called 
 a \emph{missing-common-color pair} or \emph{MCC-pair} in short with respect to the color $i$. 
 In particular, 
 when $G$ is in Conditions (a) or (b),  every vertex from $\pbar^{-1}_A(i)$
 is paired up with a vertex from $\pbar^{-1}_B(i)$. 
 
 For every MCC-pair $(a,b)$ with respect to some color $i\in[1,k]$, 
 we will exchange an alternating path $P$ from $a$ to $b$ with at most 11
 edges, where, if exist, the first, third, fifth, seventh,  ninth, and eleventh edges are uncolored
 and the second, fourth, sixth, eighth, and tenth edges are good edges colored by $i$. After $P$ is exchanged, $a$ and $b$ will be incident with edges  colored by $i$, and at most three good edges will be added to each of $R_A$ and $R_B$. With this information at hand, 
 before demonstrating the existence of such paths,  we  show that Conditions S2.1, S2.2 and S2.3  can be guaranteed  at the end of Step 2. 
 After the completion of Step 1, by~\eqref{eq-vertices-missing-i}, the total number 
 of missing colors from vertices in $A$
or from vertices in $B$ is at most $4n^{5/3}-2n$.  
 Thus there are at most $4n^{5/3}-2n$ MCC-pairs. 
 For each MCC-pair $(a,b)$ with $a,b\in V(G^*)$, 
 at most three edges will be added to each of $R_A$ and $R_B$
 when we exchange an alternating path from $a$ to $b$.  Thus 
 there will always be fewer than  
 $$3(4n^{5/3}-2n)<12 n^{5/3}$$
 edges in each of $R_A$ and $R_B$. Thus  Condition S2.1 will be satisfied at the end of Step 2. 
 And as we only ever add good edges to $R_A$ and $R_B$, Condition S2.2 will  hold automatically.
 We  now show that Condition S2.3  will also be satisfied. Recall 
 $$
 S_A=\{u\in S\cap A\,:\, d_{G_A^*}(u) \le k-2n^{2/3} \} \quad \text{and}\quad S_B=\{u\in S\cap B\,:\, d_{G_B^*}(u) \le k-2n^{2/3} \}.
 $$
 Since $\sum_{u\in A\cap S} (k-d_{G^*_A}(u)),  \sum_{u\in B\cap S} (k-d_{G^*_B}(u))\le k \le \frac{1}{2}\Delta(G)+n^{2/3}+1<2n$, it follows
 that 
 $$
 |S_A| <n^{1/3} \quad \text{and} \quad |S_B|<n^{1/3}. 
 $$
 Thus for every vertex $u\in S\setminus (S_A\cup S_B)$, $|\pbar(u)| < 2n^{2/3}$.
 For every vertex $u\in V(G^*)\setminus S$, as $d_G(u)=d_{G^*}(u) > \Delta(G)-7n^{2/3}$, it follows that $|\pbar(u)|  < k-(\frac{1}{2}(\Delta(G)-7n^{2/3})-n^{2/3})<6n^{2/3}$. 
 Thus for any $u\in V(G^*)\setminus(S_A\cup S_B)$, we have $|\pbar(u)|<6n^{2/3}$.  
 In the process of Step 2,   the number of newly colored edges of $H$ that are incident with a vertex $u\in V(G^*)\setminus (S_A\cup S_B)$ will equal  the number of alternating paths containing $u$ that have been exchanged. 
 The number of such alternating paths of which $u$ is the first vertex will equal the number of colors that missed at $u$ at the end of Step 1, which is less than $6n^{2/3}$. The number of alternating paths in which $u$ is not the first  vertex will  equal  the degree of $u$ in $R_A\cup R_B$, and so will be less than $n^{5/6}+1$. Hence the number of colored edges of $H$ that are incident with $u$ will be less than $$6n^{2/3}+n^{5/6}+1<2n^{5/6}.$$
 This applies to all vertices in $V(G^*)\setminus(S_A\cup S_B)$, and so Condition S2.3.1  will be satisfied.
 
  When $G$ is in Condition (b), for any vertex $u\in S_A\cup S_B=\{x,y\}$, since $\Delta(G)-\delta(G) \le \frac{1}{2}(1-\ve)n$,
 we have 
 $$
 |\pbar(u)|\le k-\frac{1}{2}(d_G(u)-n^{2/3}) \le k-\frac{1}{2}\left(\Delta(G)-\frac{1}{2}(1-\ve)n-n^{2/3}\right) \le \frac{1}{4}(1-\ve)n+n^{2/3}+1. 
 $$
 Hence the number of colored edges of $H$ that are incident with $u$ will be less than $$\frac{1}{4}(1-\ve)n+n^{2/3}+1+n^{5/6}+1<(\frac{1}{4}-\frac{1}{5}\ve )n.$$
 Therefore, Condition S2.3.2 will be satisfied. 
 
 When $G$ is in Condition (c), for any vertex $u\in S_A\cup S_B$, since $\Delta(G)-\delta(G) \le (1-\ve)n$,
 we have 
 $$
 |\pbar(u)|\le k-\frac{1}{2}(d_G(u)-n^{2/3}) \le k-\frac{1}{2}\left(\Delta(G)-(1-\ve)n-n^{2/3}\right) \le \frac{1}{2}(1-\ve)n+n^{2/3}+1. 
 $$
 Hence the number of colored edges of $H$ that are incident with $u$ will be less than $$\frac{1}{2}(1-\ve)n+n^{2/3}+1+n^{5/6}+1<(\frac{1}{2}-\frac{1}{3}\ve )n.$$
 Therefore, Condition S2.3.3 will be satisfied.
 
 We now show below the existence of alternating paths for MCC-pairs.  
 For a given color $i\in [1,k]$,  and vertices $a\in A$ and $b\in B$, let $N_B(a)$ be the set
 of vertices in $B$ that are joined with $a$ by an uncolored edge and are incident with a
 good edge colored $i$ such that the good edge is not incident with any vertex of $S_B$, and let $N_A(b)$ be the set of vertices in $A$ that are joined with
 $b$ by an uncolored edge and are incident with a good edge colored $i$ such that the good edge is not incident with any vertex of $S_A$. 
 In order to estimate the sizes of $N_A(b)$ and $N_B(a)$, we 
 show that $A$ and $B$ contain only a few vertices 
 that either miss the color $i$ 
 or are incident  with a non-good edge colored $i$. 
 By S2.1, there are at most  $12n^{5/3}$ edges in $R_B$, so there are fewer than $24 n^{5/6}$ vertices of
 degree at least $n^{5/6}$ in $R_B$. 
 Each non-good edge is incident with one or two 
 vertices of $R_B$ through the color $i$, 
 so there are fewer than
 $
 48 n^{5/6}
 $ 
 vertices in $B$ that are
 incident with a non-good edge colored $i$. 
 Furthermore, there are at most $2|S_B| \le 2n^{1/3}$ vertices in $B$
that are  either contained in $S_B$ 
 or adjacent to a vertex from $S_B$ through an edge with color $i$. 
Finally, there are fewer than $4n^{2/3}$
 vertices in $B$ that are missed by the color $i$. So the number of vertices in $B$ that are
 not incident with a good edge colored $i$  such that the good edge is not 
 incident with any vertex from $S_B$ is less than
 $$
 48 n^{5/6}+2n^{1/3}+4n^{2/3}<49 n^{5/6}. 
 $$
 By symmetry, 
 the number of vertices  in $A$ that are
 not incident with a good edge colored $i$ such that the good edge is not 
 incident with any vertex from $S_A$ is less than
 $49 n^{5/6}$. By S2.3.1, 
 when $\{a,b\}\cap (S_A\cup S_B) = \emptyset$, 
 \begin{equation}\tag{S2.I}\label{eqn4}
 |N_A(b)|, |N_B(a)| \ge \frac{1}{2}\left((1+4\ve/5)n-n^{2/3}\right)-2n^{5/6}-49n^{5/6}>(\frac{1}{2}+\frac{1}{3}\ve)n. 
 \end{equation} 
 When $G$ is in Condition (b) and $\{a,b\}\cap \{x,y\} \ne \emptyset$, by S2.3.2, we have 
 \begin{equation}\tag{S2.II}\label{eqn6}
 |N_A(b)|, |N_B(a)| \ge \frac{1}{2}\left((1/2+3\ve/2)n-n^{2/3} \right)-(\frac{1}{4}-\frac{1}{5}\ve)n-49n^{5/6}>\frac{3}{4}\ve n. 
 \end{equation}
 When $G$ is in Condition (c) and $\{a,b\}\cap (S_A\cup S_B) \ne \emptyset$, by S2.3.3, we have 
 \begin{equation}\tag{S2.III}\label{eqn5}
 |N_A(b)|, |N_B(a)| \ge \frac{1}{2}\left((1+\ve)n-n^{2/3}\right)-(\frac{1}{2}-\frac{1}{3}\ve)n-49n^{5/6}>\frac{1}{2}\ve n. 
 \end{equation}
 
 Let $M_B(a)$ be the set of vertices in $B$ that are joined with a vertex in $N_B(a)$ by an edge
 of color $i$, and let $M_A(b)$ be the set of vertices in $A$ that are joined with a vertex in $N_A(b)$
 by an edge of color $i$.  Note that $(S_A\cup S_B)\cap ( M_A(b) \cup M_B(a))=\emptyset$ by the choice of $N_A(b)$ and $N_B(a)$. 
 Note also that $|M_B(a)|=|N_B(a)|$     
 but some vertices
 may be in both. Similarly
 $|M_A(b)|=|N_A(b)|$. 
 
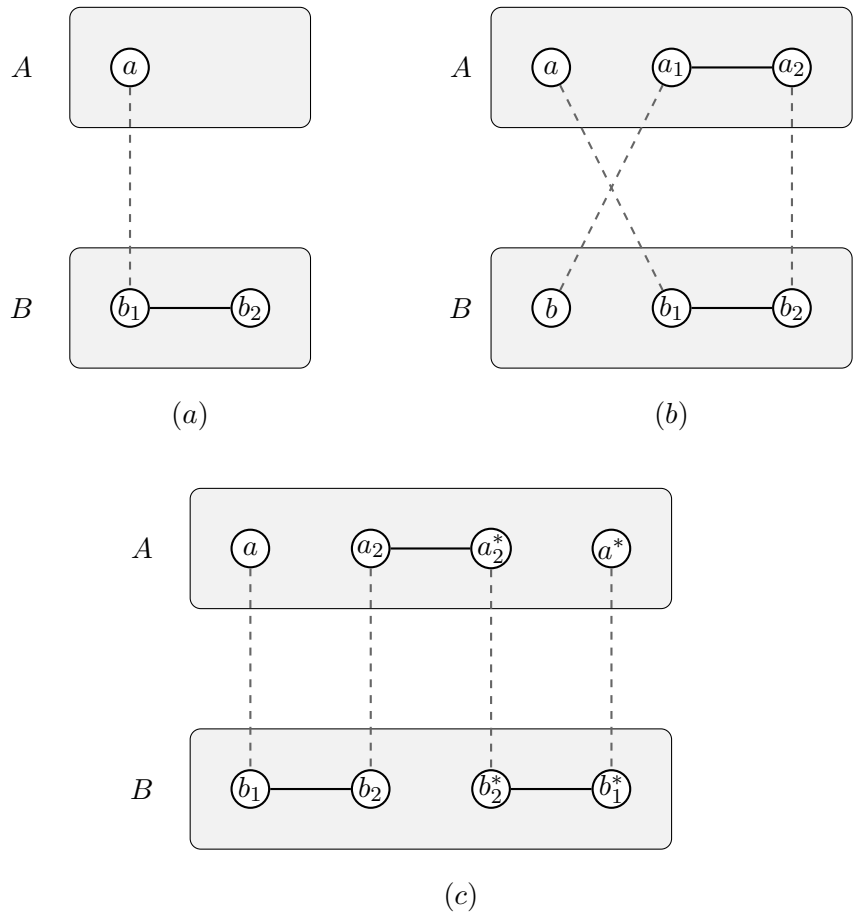
\begin{figure}[!htb]
	\begin{center}
		
		\begin{tikzpicture}[scale=0.8]
		
			\begin{scope}[shift={(-10,0)}]
		\draw[rounded corners, fill=white!90!gray] (8, 0) rectangle (12, 2) {};
		
		\draw[rounded corners, fill=white!90!gray] (8, -4) rectangle (12, -2) {};
		
		{\tikzstyle{every node}=[draw ,circle,fill=white, minimum size=0.5cm,
			inner sep=0pt]
			\draw[black,thick](9,1) node (c)  {$a$};
			\draw[black,thick](9,-3) node (d)  {$b_{1}$};
			\draw[black,thick](11,-3) node (d1)  {$b_{2}$};
			
		}
		\path[draw,thick,black!60!white,dashed]
		(c) edge node[name=la,pos=0.7, above] {\color{blue} } (d)
		;
		
		\path[draw,thick,black]
		(d) edge node[name=la,pos=0.7, above] {\color{blue} } (d1)
		;
		\node at (7.2,1) {$A$};
		\node at (7.2,-3) {$B$};
		\node at (10,-4.8) {$(a)$};	
		\end{scope}	
		
		\begin{scope}[shift={(5,0)}]
		\draw[rounded corners, fill=white!90!gray] (0, 0) rectangle (6, 2) {};
		
		\draw[rounded corners, fill=white!90!gray] (0, -4) rectangle (6, -2) {};
		
		{\tikzstyle{every node}=[draw ,circle,fill=white, minimum size=0.5cm,
			inner sep=0pt]
			\draw[black,thick](1,1) node (a)  {$a$};
			\draw[black,thick](3,1) node (a1)  {$a_1$};
			\draw[black,thick](5,1) node (a2)  {$a_2$};
			\draw[black,thick](1,-3) node (b)  {$b$};
			\draw[black,thick](3,-3) node (b1)  {$b_1$};
			\draw[black,thick](5,-3) node (b2)  {$b_2$};
		}

		\path[draw,thick,black!60!white,dashed]
		(a) edge node[name=la,pos=0.7, above] {\color{blue} } (b1)
		(a2) edge node[name=la,pos=0.7, above] {\color{blue} } (b2)
		(b) edge node[name=la,pos=0.6,above] {\color{blue}  } (a1)
		;
		
		\path[draw,thick,black]
		(a1) edge node[name=la,pos=0.7, above] {\color{blue} } (a2)
		(b1) edge node[name=la,pos=0.7, above] {\color{blue} } (b2)
		;
		
		\node at (-0.5,1) {$A$};
		\node at (-0.5,-3) {$B$};
		\node at (3,-4.8) {$(b)$};

		\end{scope}

		\begin{scope}[shift={(-8,-8)}]
		\draw[rounded corners, fill=white!90!gray] (8, 0) rectangle (16, 2) {};
		
		\draw[rounded corners, fill=white!90!gray] (8, -4) rectangle (16, -2) {};
		
		{\tikzstyle{every node}=[draw ,circle,fill=white, minimum size=0.5cm,
			inner sep=0pt]
			\draw[black,thick](9,1) node (c)  {$a$};
			\draw[black,thick](11,1) node (c1)  {$a_{2}$};
			\draw[black,thick](13,1) node (c2)  {$a^*_{2}$};
			\draw[black,thick](15,1) node (c3)  {$a^*$};
			\draw[black,thick](9,-3) node (d)  {$b_{1}$};
			\draw[black,thick](11,-3) node (d1)  {$b_{2}$};
			\draw[black,thick](13,-3) node (d2)  {$b^*_{2}$};
			\draw[black,thick](15,-3) node (d3)  {$b^*_{1}$};
			
		}
		\path[draw,thick,black!60!white,dashed]
		(c) edge node[name=la,pos=0.7, above] {\color{blue} } (d)
		(c1) edge node[name=la,pos=0.7, above] {\color{blue} } (d1)	
		(c3) edge node[name=la,pos=0.7, above] {\color{blue} } (d3)	
		(c2) edge node[name=la,pos=0.7, above] {\color{blue} } (d2)	
		;
		
		\path[draw,thick,black]
		(d) edge node[name=la,pos=0.7, above] {\color{blue} } (d1)
		(c2) edge node[name=la,pos=0.7, above] {\color{blue} } (c1)
		(d3) edge node[name=la,pos=0.7, above] {\color{blue} } (d2)
		;
		\node at (7.2,1) {$A$};
		\node at (7.2,-3) {$B$};
		\node at (12.5,-4.8) {$(c)$};	
		\end{scope}	
		
		\end{tikzpicture}
		-	  	\end{center}
	\caption{The alternating path $P$. Dashed lines indicate uncoloured edges, and solid
		lines indicate edges with color $i$.}
	\label{f1}
\end{figure}

 For a MCC-pair $(a,b)$, in order to have a uniform discussion as in the case that $\{a,b\}\cap (S_A\cup S_B) = \emptyset$,  if necessary, by exchanging an alternating path of length 2 from $a$ to another vertex $a^*$,
 and exchanging an alternating path from $b$ to another vertex $b^*$, 
 we will replace the pair $(a,b)$ by $(a^*,b^*)$ such that $\{a^*,b^*\}\cap (S_A\cup S_B)=\emptyset$.  
 Precisely, we will implement the following 
 operations to vertices in $S_A\cup S_B$. 
 For any vertex $a\in S_A$, and for each color $i\in \pbar(a)$, we take an edge $b_1b_2$ with $b_1\in N_B(a)$ and $b_2\in M_B(a)$ such that $b_1b_2$ is colored by $i$, where the edge $b_1b_2$ exists by~\eqref{eqn6}-\eqref{eqn5} and the fact that $|M_B(a)|=|N_B(a)|$.   Then we exchange the path $ab_1b_2$ by coloring $ab_1$
 with $i$ and uncoloring the edge $b_1b_2$ (See Figure~\ref{f1} (a)).
  After this, the edge $ab_1$  of $H$ is now colored by $i$, and the uncolored edge $b_1b_2$ is added to $R_B$.
 We then update the original MCC-pair that contains $a$ with respect to the color $i$ by replacing the vertex $a$ with $b_2$. 
 We do this at the vertex $a$ for every color $i\in \pbar(a)$ and then repeat the same process for every vertex in $S_A$. 
 Similarly, 
 for any vertex $b\in S_B$, and for each color $i\in \pbar(b)$, we take an edge $a_1a_2$ with $a_1\in N_A(b)$ and $a_2\in M_A(b)$ such that $a_1a_2$ is colored by $i$, where the edge $a_1a_2$ exists by~\eqref{eqn6}-\eqref{eqn5} and the fact that $|M_B(b)|=|N_B(b)|$.   Then we exchange the path $ba_1a_2$ by coloring $ba_1$
 with $i$ and uncoloring the edge $a_1a_2$. The same, we update the original MCC-pair that contains $b$ with respect to the color $i$ by replacing the vertex $b$ with $a_2$.

After the procedure above,  we have now three types MCC-pair $(u,v)$: $u,v\in A$, $u,v\in B$, and $A$
contains exactly one of $u$ and $v$ and $B$ contains the other. However, in either case, $\{u,v\}\cap (S_A\cup S_B)=\emptyset$. We will exchange alternating path for each of such pairs. 

We deal with each of the colors from $[1,k]$ in turn. 
Let $i\in [1,k]$ be a color.  We consider first an MCC-pair $(a,a^*)$
with respect to  $i$ such that $a,a^* \in A$. 
By~\eqref{eqn4}, we have  $|M_B(a^*)|>(\frac{1}{2}+\frac{1}{3}\ve)n$.
We take an edge $b_1^*b_2^*$ colored by $i$ with $b_1^* \in N_B(a^*)$ and  $b_2^*\in M_B(a^*)$.  
Then again, by~\eqref{eqn4}, we have 
$
 |M_B(a)|,  |M_A(b_2^*)|>(\frac{1}{2}+\frac{1}{3}\ve)n.  
$
Therefore, as each vertex  $c\in M_A(b_2^*)$ 
satisfies $|N_B(c)|>(\frac{1}{2}+\frac{1}{3}\ve)n$,  we have $|N_B(c)\cap M_B(a)| \ge \frac{2}{3}\ve n$. We take $a_2a_2^*$ colored by $i$ with $ a_2^*\in N_A(b_2^*)$ and  $a_2\in   M_A(b_2^*)$. 
Then we let $b_2\in N_B(a_2)\cap M_B(a)$, and let $b_1$
be the vertex in $N_B(a)$ such that $b_1b_2$ is colored by $i$. 
Now we get the alternating path $P=ab_1b_2 a_2 a_2^* b_2^* b_1^* a^*$ (See Figure~\ref{f1} (c)).  
We exchange $P$ by coloring $ab_1, b_2a_2, a_2^*b_2^*$ and $b_1^*a^*$
with color $i$ and uncoloring the edges $b_1b_2, b_1^*b_2^*$ and $a_2a_2^*$. 
After the exchange,  the color $i$ appears on edges incident with $a$ and $a^*$,
the edges $b_1b_2$ and $b_1^*b_2^*$ are added to $R_B$
and the edge $a_2a_2^*$ is added to $R_A$. We added at most 
one edge to each of $R_A$ and $R_B$ when we updated the original MCC-pair corresponding to $(a,a^*)$. Thus  we added at most three edges to each of $R_A$
and $R_B$ when we modify $\varphi$ to have the color $i$ present at both of the  vertices in the original MCC-pair corresponding to $(a,a^*)$. 
By symmetry, we can  deal with an MCC-pair $(b,b^*)$
with respect to  $i$ such that $b,b^* \in B$ similarly as above. 

Thus we consider an MCC-pair $(a,b)$
with respect to  $i$ such that $a\in A$ and $b\in B$. 
By~\eqref{eqn4}, we have 
$
|M_B(a)|,  |M_A(b)|>(\frac{1}{2}+\frac{1}{3}\ve)n.  
$
We choose $a_1a_2$ with color $i$ such that $a_1\in N_A(b)$
and $a_2\in M_A(b)$. Now as $|M_B(a)|, |N_B(a_2)| >(\frac{1}{2}+\frac{1}{3}\ve)n$ by~\eqref{eqn4}, we know that $N_B(a_2)\cap M_B(a) \ne \emptyset$. We choose $b_2\in N_B(a_2)\cap M_B(a)$
and let $b_1\in N_B(a)$ such that $b_1b_2$ is colored by $i$. 
Then $P=ab_1b_2a_2a_1 b$ is an alternating path from $a$ to $b$ (See Figure~\ref{f1} (b)). We exchange $P$ by coloring $ab_1, b_2a_2$ and $a_1b$
with color $i$ and uncoloring the edges $a_1a_2$ and $b_1b_2$. 
After the exchange,  the color $i$ appears on edges incident with $a$ and $b$,
the edge $a_1a_2$ is added to $R_A$
and the edge $b_1b_2$ is added to $R_B$. We added at most 
one edge to each of $R_A$ and $R_B$ when we updated the original MCC-pair corresponding to $(a,b)$. Thus  we added at most three edges to each of $R_A$
and $R_B$ when we modify $\varphi$ to have the color $i$ present at both of the vertices in the original MCC-pair corresponding to $(a,b)$. 
 By finding
such paths  for all MCC-pairs with respect to the color $i$,  we can increase the number of
edges colored $i$ until the color class is a 1-factor of $G^*$. By doing this for all colors,
we can make each of the $k$ color classes a 1-factor of $G^*$.

\begin{center}
	Step 3: Coloring $R_A$ and $R_B$ and extending the new color classes
\end{center}

Each of the color classes for the colors from $[1,k]$ is now a 1-factor of $G^*$. We 
now consider the multigraphs $R_A$ and $R_B$ that consist of the uncolored edges of $G^*_A$ and $G^*_B$. 
By Condition S2.1, $R_A$ and $R_B$ each has fewer than $12 n^{5/3}$ edges, and $\Delta(R_A), \Delta(R_B)< n^{5/6}+1$.  Note that $R_A$
and $R_B$ might contain parallel edges 
with endvertices in $S$. 
 By Theorem~\ref{chromatic-index} and  Theorem~\ref{lem:equa-edge-coloring},  $R_A$ and $R_B$ each have an
equalized edge-coloring with exactly $\ell:= \lceil 2n^{5/6}  \rceil$ colors $k+1, \ldots, k+\ell$. 

If $G$ is in Conditions (a) or (b), then we have $e(R_A)=e(R_B)$. 
Under these two conditions, by renaming some color classes of $R_A$ if necessary,  we can assume that in the edge colorings of $R_A$ and $R_B$, each color
appears on the same number of edges in $R_A$ as it does in $R_B$. 
When $G$ is in Condition (c), by our assumption that $G^*_B$
has more edges than $G^*_A$ does, we have  $e(R_A) \le e(R_B)$. 
In this case, we can assume that in the edge colorings of $R_A$ and $R_B$, the number of edges with a color $i\in [k+1,k+\ell]$ in $R_B$
is at least the number of edges with a color $i\in [k+1,k+\ell]$ in $R_A$.

There are fewer than $12 n^{5/3}$ edges in each of $R_A$ and $R_B$, and $\ell>n^{5/6}$, so each of
the color $i\in [k+1,k+\ell]$ appears on fewer than
$12 n^{5/6} +1$
edges in each of $R_A$ and $R_B$. We will now color some of the edges of $H$ with the $\ell$
colors from $[k+1,k+\ell]$ so that each of these color classes present at 
vertices from $V(G^*)\setminus V_\delta$. 
We  perform the following procedure for each of the $\ell$ colors in turn.

 Given a
color $i$ with $i\in [k+1,k+\ell]$, we let $A_i$ and $B_i$ be the sets of vertices in $A$ and
$B$ respectively that are incident with edges colored $i$. 
Note that $|A_i| \le  |B_i| <2(12 n^{5/6} +1) $ 
 as $R_A$ and $R_B$ each contains fewer than $ 12 n^{5/6} +1$ edges colored $i$.   Note that if $G$ is in Conditions (a) or (b), we have $|A_i|=|B_i|$;
 and we might have $|B_i| \ge |A_i|$ when $G$ is in Condition (c). 
 When $G$ is in Condition (c) and $|B_i| > |A_i|$, we let 
 $$
 A_i^* \subseteq (V_\delta\cap A) \setminus A_i 
 $$
 such that $|A_i^*|+|A_i|=|B_i|$, and just let $A_i^*=\emptyset$ otherwise.  
 Note that such $A_i^*$ exists as $|V_\delta\cap A| \ge \frac{1}{2} n^{6/7}-1$
 and $|A_i|, |B_i| <2(12 n^{5/6} +1)$. 
Let $H_i$ be the subgraph of $H$ obtained by
deleting the vertex sets $A_i \cup A_i^*$ and $B_i$ and removing all colored edges. We will show next that $H_i$ has a perfect matching and we will color 
the edges in the matching by the color $i$.

Each vertex in $V(G^*) \setminus (S_A\cup S_B)$ is incident with fewer than
$2n^{5/6}+\ell \le 5 n^{5/6}$ 
edges of $H$ that are coloured, since fewer than $2n^{5/6}$ were colored in Step 2  by S2.3.1 and at
most $2n^{5/6}+2<3n^{5/6}$ have been colored in Step 3. Also 
each vertex in $G^*$ has fewer than $ 2(12 n^{5/6} +1)$ edges that join it with a vertex in $A_i$ or $B_i$. So each vertex from $V(H_i)\setminus (S_A\cup S_B)$
is adjacent in $H_i$ to  more than 
$$
 \frac{1}{2}\left((1+4\ve/5)n-n^{2/3}\right)-5n^{5/6}-2(12 n^{5/6} +1)>\frac{1}{2}(1+\ve/2)n 
$$
vertices.

When $G$ is in condition (b), each vertex in $S_A\cup S_B$ is incident with fewer than
$(\frac{1}{4}-\frac{1}{5}\ve )n+3n^{5/6}$ 
edges of $H$ that are coloured, since fewer than $(\frac{1}{4}-\frac{1}{5}\ve )n$ were colored in Step 2 by S2.3.2 and at
most $3n^{5/6}$ have been colored in Step 3. Also 
each vertex in $G^*$ has fewer than $ 2(12 n^{5/6} +1)$ edges that join it with a vertex in $A_i$ or $B_i$. So  when $G$ is in Condition (b), each vertex from $S_A\cup S_B$
is adjacent in $H_i$ to  more than 
$$
\frac{1}{2}\left((1/2+3\ve/2)n-n^{2/3} \right)-\left((\frac{1}{4}-\frac{1}{5}\ve )n+3n^{5/6} \right)-2(12 n^{5/6} +1)> \frac{3}{4}\ve n 
$$
vertices.

When $G$ is in condition (c), each vertex in $S_A\cup S_B$ is incident with fewer than
$(\frac{1}{2}-\frac{1}{3}\ve )n+3n^{5/6}$ 
edges of $H$ that are coloured, since fewer than $(\frac{1}{2}-\frac{1}{3}\ve )n$ were colored in Step 2  by S2.3.3 and at
most $3n^{5/6}$ have been colored in Step 3. Also 
each vertex in $G^*$ has fewer than $ 2(12 n^{5/6} +1)$ edges that join it with a vertex in $A_i\cup A_i^*$ or $B_i$.  
So  when $G$ is in Condition (c), each vertex from $S_A\cup S_B$
is adjacent in $H_i$ to  more than 
$$
\frac{1}{2}\left((1+\ve)n-n^{2/3}\right)-\left((\frac{1}{2}-\frac{1}{3}\ve )n+3n^{5/6} \right)-2(12 n^{5/6} +1)> \frac{1}{2}\ve n 
$$
vertices.

Thus $\delta(H_i) \ge \frac{1}{2}\ve n$ in either case and 
$H_i$ has at most $|S_A\cup S_B| \le 2n^{1/3} <\frac{1}{2} \ve n$
vertices of degree less than $\frac{1}{2}n$. 
So $H_i$
has a 1-factor $F$ by Lemma~\ref{lem:matching-in-bipartite}. If we color the edges of $F$ with the color
$i$, then every vertex in $V(G^*)\setminus A_i^*$ is incident with an edge of color $i$. 
We repeat this procedure for each of the colors from $[k+1,k+\ell]$. After this has been
done, each of these $\ell$ colors presents at all vertices from $V(G^*)\setminus V_\delta$. So at the conclusion of Step 3,
all of the edges in $G^*_A$ and $G^*_B$ are colored, some of the edges of $H$ are colored, 
each of the $k$ color classes  for colors from $[1,k]$ is a 1-factor of $G^*$,
and  each of the $\ell$ colors  from $[k+1,k+\ell]$ 
presents at all vertices from $V(G^*)\setminus V_\delta$.

\begin{center}
	Step 4: Coloring the graph $R$ 
\end{center}
Let $R$ be the subgraph of $G^*$ consisting of the remaining uncolored edges. These
edges all belong to $H$, so $R$ is a subgraph of $H$ and hence is bipartite. 
We claim that $\Delta(R) = \Delta(G^*)-k-\ell$. 
Note that every vertex from $V(G^*)\setminus V_\delta$ 
presents every color from $[1,k+\ell]$
and so those vertices have degree  exactly  $\Delta(G^*)-k-\ell$ 
in $R$. 
For the vertices from $V_\delta$, they present all
the colors from $[1,k]$. Thus by~\eqref{eqn2}, those vertices 
have degree at most 
$$
\Delta(G^*)-\frac{1}{3} n^{6/7} -k <\Delta(G^*)-k-\ell 
$$
in $R$. 
By Theorem~\ref{konig} we can color the
edges of $R$ with $\Delta(R)$ colors from $[k+\ell+1, \Delta(G^*)]$. 
Thus $\chi'(G^*) \le k+\ell +(\Delta(G^*)-k-\ell)=\Delta(G^*)$
and so $\chi'(G^*)=\Delta(G^*)$, as desired. 

Lastly, we check that there is a polynomial time algorithm to obtain an
edge $\Delta(G)$-coloring of $G$.
By Lemma~\ref{lem:partition}, we can obtain a desired partition  
$\{A,B\}$ of $V(G)$ in polynomial time. Also,  it is polynomial time 
to edge color $G_{A}$ and $G_B$ 
by an algorithm described in~\cite{MR1156837}. Modifying $G_A$
and $G_B$ into $G_A^*$ and $G_B^*$ and the corresponding edge colorings
into equalized edge-colorings can be done in polynomial time too.  
In Step 2, the construction of the alternating paths  and swaps of the colors on the paths can be done in $O(n^3)$-time, as the total number of colors missing at vertices is $O(n^2)$ and  it takes  $O(n)$-time  to find an alternating path for a MCC-pair.  
In Step 3, there is polynomial time algorithm (see e.g. \cite{MR875324}) to edge color $R_A$ and $R_B$ using at most $\ell$ colors. Then by 
doing Kempe changes as  mentioned in the comments immediately after Theorem~\ref{lem:equa-edge-coloring}, these edge colorings can be modified
into  equalized edge-colorings in polynomial time. 
The last step is to edge color the bipartite graph $R$ using $\Delta(R)$ colors, which can be done in polynomial-time in $n$, for example, using an algorithm from~\cite{MR664720}. 
Thus, there is a polynomial time 
algorithm that gives an edge coloring of  $G$ using $\Delta(G)$ colors. 
\qed

\section*{Acknowledgment}
The second author wishes to thank Dr. Mark Ellingham 
for bringing up  reference~\cite{MR0411988} 
to her attention.



\end{document}